\documentclass[a4paper,12pt]{article}
\usepackage{mathrsfs}
\usepackage{amsmath,amsthm,amssymb}
\usepackage[dvips]{graphicx}
\makeatletter
\@addtoreset{equation}{section}

\makeatother
\title{\bf The order of complex numbers}
\author{Sun Daochun, Gu Zhendong, Liu Weiqun, Yue Chao}
\date{}
\begin{document}
\maketitle
\begin{ }

 {\bf Abstract: }In this paper, we define an ordering relation for a set of complex numbers, and research the properties and theorems
of the ordering, solve some simple complex inequalities with the
ordering.

{\bf Key words: }complex number, ordering, inequality.
\end{ }

\section{Introduction}
\mbox{}\indent In general, it is impossible to judge the order of
two complex numbers. In fact, how to judge the order of two complex
numbers is how to define  their sequencing. We think it is possible
to define the order of two complex numbers. Since the order of
complex numbers can't bring more useful interest than that of real
numbers, no one is willing to define the order of the complex
numbers.

It is assumed that the properties of the inequalities of real
numbers hold.  In this paper, we define the order of the complex
numbers, which is the extension of the order of real number, then
get some properties and theorem about this order in complex. At the
 end, using this order, we resolve some inequalities.

Let small letters denote real numbers, and let capital letters
denote complex numbers. We use $\mathscr{R}$ and $\mathscr{C}$ to
denote the set of real numbers and the set of complex numbers
respectively. We also use Greece small letter to denote the argument
of complex number, and use script capital letter to denote the
number set.

\section{The ordering relation ($\leq $) of the set of complex numbers}
\mbox{}\indent Let complex numbers $A=a_1+a_2i, B=b_1+b_2i,
C=c_1+c_2i, \cdots$

{\bf Definition 1} \quad {\it $A\leq B\Leftrightarrow(a_1< b_1)\cup
[(a_1=b_1)\cap (a_2\leq b_2)]$, that is $a_1<b_1$, or $a_1=b_1$ and
$a_2\leq b_2$.

$A\geq B\Leftrightarrow(a_1> b_1)\cup [(a_1=b_1)\cap (a_2\geq
b_2)]$. }

From the above definition, we get the following properties.

{\bf Property 1} \quad {\it(reflexivity)£º $A\leq A$. }

{\bf Property 2} \quad {\it $[(A\geq B) ~~and ~~(A\leq B)]
\Leftrightarrow$ ($A=B$). }

{\bf Proof} \quad $\{A\geq B\}\cap \{A\leq B\}$
$$ {\Leftrightarrow}
 \{(a_1> b_1)\cup [(a_1=b_1)\cap (a_2\geq b_2)]\}\cap \{(a_1< b_1)\cup [(a_1=b_1)\cap (a_2\leq b_2)]\}$$
$${\Leftrightarrow}  \{(a_1> b_1)\cap (a_1< b_1)\}\cup
 \{[(a_1=b_1)\cap (a_2\geq b_2)]\}\cap  [(a_1=b_1)\cap (a_2\leq b_2)]\}\hspace{2cm}$$
$$\hspace{2cm}\cup  \{(a_1> b_1)\cap  [(a_1=b_1)\cap (a_2\leq b_2)]\}\cup  \{[(a_1=b_1)\cap (a_2\geq b_2)]\}\cap (a_1< b_1)\}$$
$${\Leftrightarrow}  \{[(a_1=b_1)\cap (a_2\geq b_2)]\}\cap  [(a_1=b_1)\cap (a_2\leq b_2)]\}$$
$${\Leftrightarrow}  (a_1=b_1)\cap (a_2=b_2){\Leftrightarrow}  (A=B).$$

{\bf Property 3} \quad {\it(transitivity)£º $A\leq B$ and $B\leq
C\Rightarrow$ $A\leq C$. }

{\bf Proof} \quad $\{A\leq B\} \cap \{B\leq C\}$
$$\Leftrightarrow
\{[a_1< b_1] \cup [(a_1=b_1)\cap (a_2\leq b_2)]\} \cap \{[b_1< c_1]
\cup [(b_1=c_1)\cap (b_2\leq c_2)]\}$$
$$\Leftrightarrow \{[a_1< b_1]\cap [b_1< c_1]\}\cup  \{[(a_1=b_1)\cap (a_2\leq b_2)] \cap [(b_1=c_1)\cap (b_2\leq c_2)]\}\hspace{5cm}$$
$$\hspace{1cm} \cup \{[a_1< b_1]\cap [(b_1=c_1)\cap (b_2\leq c_2)]\}  \cup \{[(a_1=b_1)\cap (a_2\leq b_2)] \cap [b_1< c_1]\}
$$
$$\Rightarrow \{[a_1< c_1]\}
 \cup \{[a_1< c_1)\cap (b_2\leq c_2)]\} \cup \{[(a_1< c_1)\cap (a_2\leq b_2)] \}
  \cup  \{(a_1=c_1)\cap (a_2\leq c_2) \}$$
$$\Rightarrow \{[a_1< c_1]\} \cup  \{(a_1=c_1)\cap (a_2\leq c_2) \}\Leftrightarrow \hspace{0.5cm}A\leq C.$$

For any two complex numbers $A,B\in \mathscr{C}$, at least one
between the two relation expressions $A\leq B$, $A\geq B$ holds.
Combing with properties 1,2,3, we know that the set of complex
numbers $\mathscr{C}$ is a totally ordered set under the condition
of the above defined ordering relation $\leq$. It is easy to see
that the above defined ordering relation $\leq$ is the extension of
that of real numbers.

{\bf Theorem 1} \quad {\it $A\leq B$ $\Leftrightarrow$ $A+C\leq
B+C$.}

{\bf Proof} \quad $(A+C\leq B+C)$
$$\stackrel{D.1}{\Leftrightarrow}\hspace{0.3cm}\{(a_1+c_1< b_1+c_1) \cup [(a_1+c_1=b_1+c_1)\cap (a_2+c_2\leq b_2+c_2)]\}$$
$$\stackrel{R.P}{\Leftrightarrow}\hspace{0.3cm}\{(a_1< b_1) \cup [a_1=b_1)\cap (a_2\leq b_2)]\}\hspace{0.3cm}
\stackrel{D.1}{\Leftrightarrow}\hspace{0.3cm}(A\leq B).$$ The
notation $D.1$ in the first equivalent symbol denote the following
result is deduced by definition 1. The notation $R.P$ in the second
equivalent symbol denote the following result is deduced by the
properties of the inequalities of real numbers. In the following, we
also use the similar remark.

{\bf Theorem 2} \quad {\it The term in the inequalities can be
moved from one side to the other, that is $A+B\leq
C\hspace{0.3cm}\Leftrightarrow\hspace{0.3cm}A\leq C-B.$ }

{\bf Proof} \quad $A+B\leq
C\hspace{0.3cm}\stackrel{T.1}{\Leftrightarrow}\hspace{0.3cm}A+B-B\leq
C-B \hspace{0.3cm}\Leftrightarrow\hspace{0.3cm}A\leq C-B.$

{\bf Theorem 3} \quad {\it $A\leq B$, and $C\leq D\Rightarrow$
$A+C\leq B+D$.}

{\bf Proof} \quad Since $A\leq B$
$\stackrel{T.1}{\Leftrightarrow}\hspace{0.3cm}$ $A+C\leq B+C$ and
$C\leq D$ $\stackrel{T.1}{\Leftrightarrow}\hspace{0.3cm}$ $B+C\leq
B+D$, by the transitivity, we get $A+C\leq B+D$.

{\bf Theorem 4} \quad {\it Let $r>0$, then $A\leq B$
$\Leftrightarrow$ $rA\leq rB$. }

{\bf Proof} \quad
$(rA\leq rB)\stackrel{D.1}{\Leftrightarrow}\hspace{0.3cm}(ra_1< rb_1)\cup [(ra_1=rb_1)\cap (ra_2\leq rb_2)]$\\
$\stackrel{R.P}{\Leftrightarrow}\hspace{0.3cm}(a_1< b_1)\cup
[(a_1=b_1)\cap (a_2\leq b_2)]$
$\stackrel{D.1}{\Leftrightarrow}\hspace{0.3cm}$$(A\leq B)$.

{\bf Theorem 5} \quad {\it Let $r<0$, then $(A\leq B)
\Leftrightarrow rA\geq rB$. Especially, $(A\leq B)\Leftrightarrow
(-A\geq -B)$. }

{\bf Proof} \quad $(rA\geq rB) \stackrel{D.1}{\Leftrightarrow}\hspace{0.3cm}(ra_1>rb_1)\cup [(ra_1=rb_1)\cap (ra_2\geq rb_2)]$\\
 $\stackrel{R.P}{\Leftrightarrow}\hspace{0.3cm}$ $(a_1< b_1)\cup [(a_1=b_1)\cap (a_2\leq b_2)]$
$\stackrel{D.1}{\Leftrightarrow}\hspace{0.3cm}$ $(A\leq B)$.

\section{The operation of the set}
\mbox{}\indent {\bf Definition 2} \quad {\it Let the set of complex
numbers $\mathscr{B}\subset \mathscr{C}$, $\theta\in \mathscr{R}$.
We define the rotation set of $\mathscr{B}$ by
$e^{i\theta}\mathscr{B}:=\{Ze^{i\theta}\in \mathscr{C}; Z\in
\mathscr{B}\}$. }

{\bf Remark 1} \quad $e^{i\theta}\mathscr{B}$ denote the set
$\mathscr{B}$ rotates $\theta$ radian around the origin.

{\bf Theorem 6} \quad {\it $We^{i\theta}\in \mathscr{B}$
$\Leftrightarrow$ $W\in e^{-i\theta}\mathscr{B}$. }

{\bf Proof} \quad Let $W:=Ze^{-i\theta}$ $\Leftrightarrow$
$Z=We^{i\theta}$. By the definition, we get
$$e^{-i\theta}\mathscr{B}:=\{Ze^{-i\theta}\in \mathscr{C}; Z\in \mathscr{B}\}
=\{W\in \mathscr{C}; We^{i\theta}\in \mathscr{B}\}
$$
From the last equality we deduce that $W\in e^{-i\theta}\mathscr{B}$
$\Leftrightarrow$ $We^{i\theta}\in \mathscr{B}$.

{\bf Example 1} \quad When $\mathscr{B}=\{Z; Z\geq A\in
\mathscr{C}\}$ is a semi-open and semi-closed half plane,
$e^{i\theta}\mathscr{B}$ is also a semi-open and semi-closed half
plane, that is ($A=a_1+a_2i,Z=z_1+z_2i$)
$$e^{i\theta}\mathscr{B}\stackrel{D.2}{=}\{Ze^{i\theta}; Z\geq A\}=\{Z; e^{-i\theta}Z\geq A\}$$
$$=\{Z; (z_1\cos \theta+z_2\sin \theta)+(z_2\cos \theta-z_1\sin \theta)i\geq a_1+a_2i\}$$
$$\stackrel{D.1}{=}(Z; z_1\cos \theta+z_2\sin \theta > a_1)\cup
[(z_1\cos \theta+z_2\sin \theta=a_1)\cap (z_2\cos \theta-z_1\sin
\theta\geq a_2)].$$ The first part is an open half plane, the part
in $[ . ]$ is a half line.

{\bf Definition 3} \quad {\it Let the set of complex numbers
$\mathscr{B}\subset \mathscr{C}$, $r\geq 0$. We define the dilation
set of $\mathscr{B}$ by $r\mathscr{B}:=\{rZ\in \mathscr{C}; Z\in
\mathscr{B}\}$. }

{\bf Remark 2} \quad $r\mathscr{B}$ denotes $\mathscr{B}$ makes a
stretching at the extension ratio $r$.

{\bf Theorem 7} \quad {\it $rW\in \mathscr{B}$ $\Leftrightarrow$
$W\in \mathscr{B}/r$. }

{\bf Proof} \quad Let $W:=Z/r$ $\Leftrightarrow$ $Z=rW$. By the
definition, we get
$$\mathscr{B}/r:=\{Z/r\in \mathscr{C}; Z\in \mathscr{B}\}
=\{W\in \mathscr{C}; rW\in \mathscr{B}\}.$$ Namely,$W\in
\mathscr{B}/r$ $\Leftrightarrow$ $rW\in \mathscr{B}$.

{\bf Definition 4} \quad {\it Let the set of complex numbers
$\mathscr{B}\subset \mathscr{C}$, $A\in \mathscr{C}$. We define the
translation set of $\mathscr{B}$ by $\mathscr{B}+A:=\{Z+A\in
\mathscr{C}; Z\in \mathscr{B}\}$. }

{\bf Remark 3} \quad $A+\mathscr{B}$ denotes $\mathscr{B}$ is
translated.

{\bf Theorem 8} \quad {\it $A+W\in \mathscr{B}$ $\Leftrightarrow$
$W\in \mathscr{B}-A$. }

{\bf Proof} \quad Let $W:=Z-A$ $\Leftrightarrow$ $Z=W+A$. By the
definition, we get
$$\mathscr{B}-A:=\{Z-A\in \mathscr{C}; Z\in \mathscr{B}\}
=\{W\in \mathscr{C}; W+A\in \mathscr{B}\}.$$ Namely, $W\in
\mathscr{B}-A$ $\Leftrightarrow$ $W\in \mathscr{B}+A$.

{\bf Definition 5} \quad {\it Let the set of complex numbers
$\mathscr{B}\subset \mathscr{C}$. We define the inversion set of
$\mathscr{B}$ by $1/\mathscr{B}:=\{1/Z\in \mathscr{C}; Z\in
\mathscr{B}\}$. }

{\bf Remark 4} \quad $1/\mathscr{B}$ and $\mathscr{B}$ is symmetric
about the unit circumference $\{|Z|=1\}$, or is called the inversion
transform.

{\bf Theorem 9} \quad {\it $1/W\in \mathscr{B}$ $\Leftrightarrow$
$W\in 1/\mathscr{B}$. }

{\bf Proof} \quad Let $W:=1/Z$ $\Leftrightarrow$ $Z=1/W$. By the
definition, we get
$$1/\mathscr{B}:=\{1/Z\in \mathscr{C}; Z\in \mathscr{B}\}
=\{W\in \mathscr{C}; 1/W\in \mathscr{B}\}.$$ Namely, $W\in
1/\mathscr{B}$ $\Leftrightarrow$ $1/W\in \mathscr{B}$.

{\bf Example 2} \quad
 When $\mathscr{B}=\{Z; Z\geq A\in \mathscr{C}\}$ is a
semi-open and semi-closed half plane, $1/\mathscr{B}$ is a semi-open
and semi-closed disc which is symmetric about the real axis, that is
($Z:=z_1+z_2i\ne 0$)
$$1/\mathscr{B}\stackrel{D.5}{=}\{1/Z; Z\geq A\}=\{Z; 1/Z\geq A\}$$
$$=\{Z; \frac{z_1}{z^2_1+z^2_2}-i\frac{z_2}{z^2_1+z^2_2}\geq a_1+a_2i\}$$
$$\stackrel{D.1}{=}(\frac{z_1}{z^2_1+z^2_2}> a_1)\cup
[(\frac{z_1}{z^2_1+z^2_2}=a_1)\cap (-i\frac{z_2}{z^2_1+z^2_2}\geq
a_2)]$$
$$=(a_1z^2_1+a_1z^2_2> z_1)\cup
[(a_1z^2_1+a_1z^2_2=z_1)\cap (a_2z^2_1+a_2z^2_2\leq -z_2  )].$$ The
first part is an open disc, its center is $\frac{1}{2a_1}+0i$, its
radius is $\frac{1}{2a_1}$. The part in $[ . ]$ is a circular arc
whose point isn't the origin.

{\bf Definition 6} \quad {\it Let the set of complex numbers
$\mathscr{B}\subset \mathscr{C}$. We define the radication set of
$\mathscr{B}$ by $\mathscr{B}^{1/2}:=\{Z^{1/2}; Z\in \mathscr{B}\}$.
}

{\bf Theorem 10} \quad {\it $W^{2}\in \mathscr{B}$ $\Leftrightarrow$
$W\in \mathscr{B}^{1/2}$. }

{\bf Proof} \quad Let $Z^{1/2}:=\{-W,W\}$ $\Leftrightarrow$
$Z=W^{2}$. By the definition, we get
$$\mathscr{B}^{1/2}:=\{Z^{1/2};\  Z\in \mathscr{B}\}
=\{-W;\  W^{2}\in \mathscr{B}\}\cup \{W;\  W^{2}\in \mathscr{B}\}
=\{W;\ W^{2}\in \mathscr{B}\}.$$ Namely, $W\in \mathscr{B}^{1/2}$
$\Leftrightarrow$ $W^{ 2}\in \mathscr{B}$.

{\bf Example 3} \quad When $\mathscr{B}=\{Z; Z\geq A\}$ is a
semi-open and semi-closed half plane, its radication set
$\mathscr{B}^{1/2}$ is a semi-open and semi-closed domain whose
boundary is hyperbola, that is
$$\mathscr{B}^{1/2}\stackrel{D.6}{=}\{Z^{1/2}; Z\geq A\}=\{Z; Z^{2}\geq A\}$$
$$=\{Z; (z_1+z_2i)^2\geq a_1+a_2i\}$$
$$=\{Z; (z^2_1-z^2_2)+2z_1z_2i\geq a_1+a_2i\}$$
$$\stackrel{D.1}{=}
[z^2_1-z^2_2>a_1]\cup [(z^2_1-z^2_2=a_1)\cap (2z_1z_2\geq a_2)].$$
The hyperbola $z^2_1-z^2_2=a_1$ splits the plane into three parts
(When $a_1=0$, it splits the plane into four parts.). If $a_1>0$,
$\mathscr{B}^{1/2}$ is two non-neighbor parts that doesn't contain
the origin. If $a_1<0$, $\mathscr{B}^{1/2}$ is a connected part that
contains the origin. We call it the hyperbola domain.

\section{Solving inequalities}
\mbox{}\indent

{\bf Definition 7} \quad {\it To solve the inequality $f(Z)\geq
g(Z)$ means to find the set $\{Z\in \mathscr{C}; f(Z)\geq g(Z)\}$. }

{\bf Definition 8} \quad {\it Define
$\mathscr{D}(A):=\mathscr{D}(Z\geq A):=\{Z\in \mathscr{C}; Z\geq
A\}=\{Z\in \mathscr{C}; Z\in \mathscr{D}(Z\geq A)\}$.\par
$\mathscr{D}( A)$ is a semi-open and semi-closed perpendicular half
plane that is split by the perpendicular line $\{Z\in \mathscr{C};Re
Z= Re A\}$. }

{\bf 1. Solving the linear inequality with one unknown} \quad Let
$A=re^{i\theta}$($r\geq 0$), solve the inequality $AZ-B\geq 0$.

{\bf Solution}\quad  By the definition and properties of
inequalities, we get the solution set $\mathscr{S}$ is the
following.

 $$\mathscr{S}\stackrel{D.7}{=}\{Z\in \mathscr{C}; AZ-B\geq 0\}\stackrel{T.2,4}{=}\{Z\in \mathscr{C}; e^{i\theta}Z\geq B/r\}$$
$$ \stackrel{W=e^{i\theta}Z}{=}\{We^{-i\theta}\in \mathscr{C}; W\geq B/r\}
 \stackrel{D.2}{=}e^{-i\theta}\{W\in \mathscr{C}; W\geq B/r\}$$
$$ \stackrel{D.7}{=}e^{-i\theta}\mathscr{D}(Z\geq B/r).$$
$$ \stackrel{D.8}{=}e^{-i\theta}\mathscr{D}( B/r).$$

Let a perpendicular half plane rotate $\theta$ radian anticlockwise
around the origin, we get a half plane whose boundary is a oblique
line, that is the solution set $\mathscr{S}$ .\par

{\bf 2. Solving the linear inequalities with one unknown}\quad Let
$A=re^{i\theta}$,$C=ue^{i\phi}$($r,u\geq 0$),solve the inequalities
$AZ-B\geq 0$,\ $CZ-D\geq 0$.

{\bf Solution}\quad Since the solution sets of two inequalities are
$$ e^{-i\theta}\mathscr{D}( B/r\},\hspace{0.5cm}e^{-i\phi}\mathscr{D}( B/u\},$$
respectively. The solution set of the inequalities $\mathscr{S}$ is
the following.
$$\mathscr{S}= [e^{-i\theta}\mathscr{D}( B/r)]\cap [e^{-i\phi}\mathscr{D}( B/u)].$$
The solution set $\mathscr{S}$ denotes the intersection of two half
planes .

{\bf 3. Solving the linear fractional inequality}\quad Let $B-AC\ne
0$, solve the inequality $\frac{AZ+B}{Z+C}\geq D$.\par

{\bf Solution}\quad Let $B-AC:=re^{i\theta}$, then the solution set
$\mathscr{S}$ is the following.
 $$\mathscr{S}\stackrel{D.7}{=}\{Z\in \mathscr{C};\frac{AZ+B}{Z+C}\geq D\}
 \stackrel{T.1,4}{=}\{Z\in \mathscr{C};\frac{1}{e^{-i\theta}(Z+C)}\geq \frac{D-A}{r}\}$$
 $$\stackrel{D.8}{=}\{Z\in \mathscr{C};\frac{1}{e^{-i\theta}(Z+C)}\in \mathscr{D}(\frac{D-A}{r})\}
\stackrel{T.9}{=}\{Z\in \mathscr{C};e^{-i\theta}(Z+C)\in
1/\mathscr{D}(\frac{D-A}{r})\}$$
 $$\stackrel{T.6}{=}\{Z\in \mathscr{C};(Z+C)\in  \frac{1}{\mathscr{D}(\frac{D-A}{r})}\cdot e^{i\theta}\}
 \stackrel{T.8}{=}\{Z\in \mathscr{C};Z\in \frac{1}{\mathscr{D}(\frac{D-A}{r})}\cdot e^{i\theta}-C\}$$
 $${=} \frac{1}{\mathscr{D}(\frac{D-A}{r})}\cdot e^{i\theta}-C.$$

The solution set $\mathscr{S}$ is a semi-open and semi-closed disc.
Let a semi-open and semi-closed perpendicular half plane
$\mathscr{D}(\frac{D-A}{r})$ become a disc
$\frac{1}{\mathscr{D}(\frac{D-A}{r})}$ by the inversion transform
defined as example 2, then let it rotate $\theta$ radian around the
origin, we get the disc $\frac{1}{\mathscr{D}(\frac{D-A}{r})}\cdot
e^{i\theta}$, finally translating the disc we get the solution set
$\frac{1}{\mathscr{D}(\frac{D-A}{r})}\cdot e^{i\theta}-C$.

{\bf 4. Solving the inequality of the second order}\quad Let $A\ne
0$, solve the inequality of the second order $AZ^2+BZ+C\geq 0$.\par

{\bf Solution}\quad Let $A:=re^{i\theta}$, then the solution set
$\mathscr{S}$ is the following.
 $$\mathscr{S}\stackrel{D.7}{=}\{Z; AZ^2+BZ+C\geq 0\}
 \stackrel{T.1,4}{=}\{Z; e^{i\theta}(Z+\frac{B}{2A})^2\geq \frac{B^2-4AC}{4rA}\}$$
 $$\stackrel{D.8}{=}\{Z; [e^{i\theta/2}(Z+\frac{B}{2A})]^2\in \mathscr{D}(\frac{B^2-4AC}{4rA})\}
\stackrel{T.10}{=}\{Z; e^{i\theta/2}(Z+\frac{B}{2A})\in
\mathscr{D}^{1/2}(\frac{B^2-4AC}{4rA})\}$$
$$\stackrel{T.6}{=}\{Z; Z+\frac{B}{2A}\in \mathscr{D}^{1/2}(\frac{B^2-4AC}{4rA})\cdot e^{-i\theta/2}\}$$
$$\stackrel{T.8}{=}\{Z; Z\in \mathscr{D}^{1/2}(\frac{B^2-4AC}{4rA})\cdot e^{-i\theta/2}-\frac{B}{2A}\}
{=} \mathscr{D}^{1/2}(\frac{B^2-4AC}{4rA})\cdot
e^{-i\theta/2}-\frac{B}{2A}.$$

The solution set $\mathscr{S}$ is a semi-open and semi-closed
hyperbola domain. Let a semi-open and semi-closed perpendicular half
plane $\mathscr{D}(\frac{B^2-4AC}{4rA})$ become a hyperbola domain
$\mathscr{D}^{1/2}(\frac{B^2-4AC}{4rA})$ by the radication defined
as example 3, then let the hyperbola domain rotate $\theta$ radian
anticlockwise around the origin, we get the hyperbola domain
$\mathscr{D}^{1/2}(\frac{B^2-4AC}{4rA})\cdot e^{-i\theta/2}$,
finally translating it we get the solution set
$\mathscr{D}^{1/2}(\frac{B^2-4AC}{4rA})\cdot
e^{-i\theta/2}-\frac{B}{2A}$.

\

\

\

\begin{flushright}
{\it Sun Daochun}\\
School of Mathematical Sciences\\
South China Normal University\\
Guangzhou 510631\\
People's Republic of China\\
E-mail\,:\, sundch@scnu.edu.cn

\

{\it Gu Zhendong}\\
School of Mathematical Sciences\\
South China Normal University\\
Guangzhou 510631\\
People's Republic of China\\
E-mail\,:\, guzhd@qq.com

\

{\it Liu Weiqun}\\
School of Mathematic\\
Jia Ying University\\
Guangdong Meizhou 514015\\
People's Republic of China\\
E-mail\,:\, mzlwq226@21cn.com

\

{\it Yue Chao}\\
School of Mathematical Sciences\\
South China Normal University\\
Guangzhou 510631\\
People's Republic of China\\
E-mail\,:\, 1048348982@qq.com
\end{flushright}

\end{document}